\documentclass[10pt]{amsart}

\usepackage{amsfonts,amssymb,amsmath,amsthm,amscd}
\usepackage[latin2]{inputenc}
\usepackage[T1]{fontenc}
\usepackage[mathscr]{eucal}
\usepackage{indentfirst}
\usepackage{graphicx,graphics}
\usepackage{pict2e,epic}
\usepackage{MnSymbol}
\usepackage[colorlinks,citecolor=red,urlcolor=blue,bookmarks=false,hypertexnames=true]{hyperref}
\usepackage[margin=3cm]{geometry}

\allowdisplaybreaks
\numberwithin{equation}{section}

\theoremstyle{plain}
\newtheorem{Th}{Theorem}[section]

\newtheorem{Lemma}[Th]{Lemma}
\newtheorem{Cor}[Th]{Corollary}

\newtheorem*{5Theorem-non}{Proof of Proposition \ref{5proposition24}}
\newtheorem*{x5Theorem-non}{Proof of Proposition \ref{x5proposition24}}

\newtheorem*{Theorem-non2}{Theorem 4.8}

\newtheorem*{Theorem-non3}{Corollary 3.10}

\theoremstyle{definition}

\newtheorem{Rem}[Th]{Remark}

\begin{document}
\setlength{\abovedisplayskip}{3pt}
\setlength{\belowdisplayskip}{3pt}

\email{Jkim6@mcneese.edu}
\title{Approximations of $SL(3,\mathbb{Z})$ Hecke--Maass $L$-Functions by short Dirichlet polynomials}
\author{Jiseong Kim}
\address{Mcneese State University, Department of Mathematical Sciences}

\begin{abstract}
We study averages of $L$-functions associated with Hecke-Maass cusp forms for $SL(3,\mathbb{Z})$, multiplied by Dirichlet polynomials built from the Fourier coefficients of the cusp forms. To prove this, we employ a variant of the Kuznetsov trace formula. In particular, we show that the reciprocals of these $L$-functions can be approximated by very short Dirichlet polynomials, on average over $t$ and over the forms.

\end{abstract}

\maketitle
\markboth{Jiseong Kim}{Approximations of $L$-functions by short Dirichlet polynomials}

\section{Introduction}

\noindent
We consider an orthogonal basis of Maass cusp forms for $SL(3,\mathbb{Z})$ consisting of Hecke--Maass cusp forms, denoted by $\{F_j\}$. Each $F_j$ has type $v_j=(v_{j,1},v_{j,2}, v_{j,3})$, and the Langlands parameter $\mu_{j}=(\mu_{j1}, \mu_{j,2}, \mu_{j,3}):=\left(2v_{j,1}+v_{j,2}, v_{j,2}-v_{j,1},-v_{j,1}-2v_{j,2} \right). $ 
Let 
\begin{equation}\label{eq:GFj}
L_{F_j}(s)
:= \sum_{n=1}^{\infty} A_j(1,n)n^{-s}
= \prod_p \prod_{i=1}^{3}\left(1-\frac{\alpha_{j,i}(p)}{p^s}\right)^{-1}
\end{equation}
be the Godement-Jacquet $L$-function associated with $F_j$ (see \cite[Chapter~6]{Go}).For a fixed type $v_0$ satisfying
\[
|\mu_{0,i}| \asymp |v_{0,i}| \asymp U
\quad \text{for all } i,
\]
we restrict our attention to forms that are away from the Weyl chamber walls and are not self-dual forms.
Let
\[
M = U^{c}
\quad \text{for some constant } c \in (3/4,1).
\]
We use the test function $h_{U,M}(\mu)$ for $\mu=(\mu_{1},\mu_{2},\mu_{3})$ introduced in \cite{sun2024}(see also \cite{ButtBlom2020}, \cite{BZ1}), which is

\begin{align}
h_{U,M}(\mu)
&:= P(\mu)^2
\left(
\sum_{w \in \mathcal{W}}
\psi\!\left(\frac{w(\mu)-\mu_0}{M}\right)
\right)^2, \notag\\[0.3em]
\psi(\mu)
&:= \exp\!\left(-\left(\mu_1^2+\mu_2^2+\mu_3^2\right)\right), \notag\\[0.3em]
P(\mu)
&:= \prod_{0 \le n \le \mathcal{A}}
\prod_{k=1}^{3}
\frac{
\left(v_k-\tfrac{1}{3}(1+2n)\right)
\left(v_k+\tfrac{1}{3}(1+2n)\right)
}{
\lvert v_{0,k} \rvert^{2}
},
\end{align}
for sufficiently large $\mathcal{A}.$ Note that $h_{U,M}$ is nonnegative and it is essentially supported for the case when $|\mu_{j,k}| \asymp |v_{j,k}| \asymp U.$

For convenience, we define
\begin{align}
\sum^{\ast}_{F_j}
&:=\sum_{\substack{F_j}}\frac{h_{U,M}(\mu_{j})}{\left\|F_j\right\|^2 \prod_{k=1}^3 \cos \left(3 v_{j, k} / 2\right)},
\\
H:&=\sum_{F_j} \frac{h_{U,M}\left(\mu_j\right)}{\left\|F_j\right\|^2 \prod_{k=1}^3 \cos \left(3 v_{j, k} / 2\right)} \times 1= \frac{1}{192 \pi^5} \int_{\operatorname{Re}(\mu)=0} h_{U, M}(\mu) \mathrm{d}_{\mathrm{spec}} \mu
\end{align}
to simplify notation.
In \cite{sun2024}, the authors investigate an asymptotic formula for the mollified second moment of $L_{F_{j}}(s)$. They define a mollifier $M(s,F_j)$ for $L_{F_{j}}(s)$ by
\[
M(s,F_j):=\sum_{\ell \le L^2} \frac{A_j(1,\ell)}{\ell^s}\,x_{\ell},
\]
where
\[
x_{\ell}
:=\mu(\ell)\,\frac{1}{2 \pi i} \int_{(3)}
\frac{(L^2/\ell)^u-(L/\ell)^u}{u^2}\,\frac{du}{\log L}
=
\begin{cases}
\mu(\ell), & \ell \le L,\\[4pt]
\mu(\ell)\,\dfrac{\log (L^2/\ell)}{\log L}, & L<\ell \le L^2,\\[8pt]
0, & \ell>L^2.
\end{cases}
\]
Here $L=U^{\delta/2}$ for some suitably small constant $\delta>0$. 

\begin{Th}\label{squaremoment}
For sufficiently small $\delta,\theta_1>0$, we have
\[
\sum_{F_{j}}^{\ast} 
\left|L_{F_{j}}\!\left(\tfrac12+\theta+i t\right)
M\!\left(\tfrac12+\theta+i t, F_j\right)\right|^2
\le H\left(1+O\!\left(U^{-2 \theta \delta}\right)\right),
\]
uniformly for $1/\log U \le \theta \le 1/2$ and $|t|<U^{\theta_1}.$
\end{Th}
\begin{proof}See \cite[Proposition 1.4]{sun2024}
\end{proof}

Let $A,B\in\mathbb{N}$ and $\alpha\in\mathbb{N} \cup \{0\}$. Let $\tau(n)$ be a function such that $|\tau(n)| \ll (\log U)^{O(1)}$ and  $\tau(1)=1$. Define the Dirichlet polynomial
\[
M_{j;A,B}(s):=
\sum_{\substack{n_{1}<A\\ n_{2}<B}}
\frac{A_{j}(n_{2},n_{1})\,\tau(n_{1})\tau(n_{2})}{n_{1}^{s}n_{2}^{\alpha s}}.
\]
In this paper, we prove the following result.

\begin{Th}\label{Main}
Suppose that $T$ and $U$ are sufficiently large and satisfy $T^{1+2\varepsilon} \ll U$ for a fixed small $\varepsilon>0$. Suppose also that
\[
\begin{aligned}
A\,B^{\frac{3}{2}-\frac{\alpha}{2}}
&\ll U^{-1-\beta-100\varepsilon}M^{2}, \\[0.4em]
A^{\frac{39}{64}}\,B^{\frac{71-32\alpha}{64}}
&\ll U^{\,\frac{11}{64}-\beta-100\varepsilon}, \\[0.4em]
A^{\frac{1}{2}}\,B^{1-\frac{\alpha}{2}}
&\ll U^{-3/2-\beta-100\varepsilon}M^{2}, \\[0.4em]
 \max\left(1,B^{1/2-\alpha/2 - \theta  - \alpha \theta + \varepsilon}\right)
&\ll T\,U^{-\beta}, \\[0.4em]
B^{\frac{58}{42}-\frac{\alpha}{2}-\frac{2\theta}{3}-\alpha \theta}
&\ll U^{\,\frac{1}{2}+\theta-\beta},
\end{aligned}
\]

for some sufficiently small constant $\beta>0$. Then, for any positive function $\theta=\theta(U)$,
\[
\sum_{\substack{F_j}}^{\ast}
\int_{T}^{2T}
L_{F_{j}}\!\left(\tfrac12+\theta + it\right)\,
M_{j;A,B}\!\left(\tfrac12+\theta+it\right)\,dt
=
TH\left(1+O\left((\log U)^{O(1)}U^{-\beta}\right)\right).
\]
\end{Th}

\medskip
\noindent

\noindent

Note that the mollifier can be regarded as $M_{j; U^{\delta},2}(s)$ with $\alpha=0$.
Combining Theorem~\ref{squaremoment} with Theorem~\ref{Main}, we obtain the following corollary.

\begin{Cor}\label{cor:mollified} Let $\eta(U)$ be a function such that $\eta(U) \rightarrow \infty$ as $U \rightarrow \infty.$ 
For sufficiently small $\delta,\theta_1>0$ satisfying $\delta \leq \theta_{1}/2,$ we have
\[
\sum_{F_{j}}^{\ast}
\int_{T}^{2T}
\left|1-L_{F_{j}}\!\left(\tfrac12+\theta+it\right)
M\left(\tfrac12+\theta+it, F_j\right)\right|^2 dt
=
o\!\left(HT\right)
\]
uniformly for $\eta(U)/\log U \le \theta \le 1/2$ and $T=U^{\theta_1}$.
\end{Cor}

\begin{proof}
By taking $\delta \leq \theta_{1}/2$ and $\beta=\delta/2,$ it is easy to see that all the conditions in Theorem \ref{Main} are satisfied.
Expanding the square and applying Theorem~\ref{Main}, it suffices to show that
\[
\sum_{F_{j}}^{\ast}\int_{T}^{2T}
\left|L_{F_{j}}\!\left(\tfrac12+\sigma+i t\right)
M\!\left(\tfrac12+\sigma+i t, F_j\right)\right|^2 dt
\le
H\left(T+o\!\left(T\right)\right),
\]
which follows directly from Theorem~\ref{squaremoment}.
\end{proof}
\begin{Rem}
Since
\[
\frac{h_{U,M}(\mu_j)}
{\|F_j\|^{2}\prod_{k=1}^{3}\cos\!\left(\tfrac{3}{2}v_{j,k}\right)}
> 0,
\]
the corollary implies that
\[
L_{F_{j}}\!\left(\tfrac12+\sigma+it\right)^{-1}
\sim
M\!\left(\tfrac12+\sigma+it, F_j\right),
\]
on average over $t \in [T,2T]$ and over the forms $F_j$.
\end{Rem}

\noindent One may ask whether $M_{j;A,B}(s)$ can be shorter than any fixed power of $U$ and still satisfy the upper bound in Corollary~1.3. We show that this is possible, and we provide the proof in the final section.

\begin{Th}\label{last}
Let $\eta(U)$ be a function such that $(\log\log U)^2=o(\eta(U))$, and suppose that
\[
|T|=U^{\theta_1}, \qquad 
\frac{\eta(U)}{\log U}<\theta\leqslant \tfrac12,
\qquad 
A=\exp\!\left(\frac{\log U}{\eta(U)^{1/2}}\right),
\]
for some sufficiently small positive constant $\theta_1$.
Let
\[
M_{j;A,2}\!\left(\tfrac12+\theta+it\right)
=1+ \sum_{1<n_1<A} 
\frac{A_j(1,n_1)C\mu(n_1)}
{n_1^{1/2+\theta+it}}.
\] 
where 
\[
C=
\left(\frac{1-\zeta(1+2\theta)}{\zeta(1+2\theta)}
\left(\sum_{d=1}^{\infty}\frac{1}{d^{1+2\theta}}
\left(
\sum_{k=2}^{\infty}\frac{\mu(dk)}{k^{1+2\theta}}
\right)^{2}+\sum_{d=2}^{\infty}\frac{1}{d^{1+2\theta}}
\left(
\sum_{k=1}^{\infty}\frac{\mu(dk)}{k^{1+2\theta}}
\right)^{2} \right)
\right)^{-1/2}.
\]
Then 
\[
\sum_{F_{j}}^{\ast}
\int_{T}^{2T}
\left|1-L_{F_{j}}\!\left(\tfrac12+\theta+it\right)
M_{j;A,2}\left(\tfrac12+\theta+it\right)\right|^2 dt
=
o\!\left(HT\right).
\]
\end{Th}

For the remainder of the paper, we treat $\varepsilon$ as an arbitrarily small positive number without further comment. $f\asymp g$ means that $f$ and $g$ are bounded above and below by constant multiples of each other.

\medskip
\noindent

\noindent
\section{Preliminaries}

\noindent

\subsection{Hecke--Maass cusp forms over $SL(3,\mathbb{Z})$}

\noindent
The Euler product of $L_{F_j}(s)$ is
\[
L_{F_j}(s)=\prod_{p}\left(1-\frac{A_j(1,p)}{p^{s}}+\frac{A_j(p,1)}{p^{2s}}-\frac{1}{p^{3s}}\right)^{-1},
\]
and $A_j(m,n)=\overline{A_j(n,m)}$.
The $L$-function $L_{F_j}(s)$ satisfies the functional equation
\begin{equation}\label{x5functionalequation}
G_{v,j}(s)L_{F_j}(s)=\widetilde{G}_{v,j}(1-s)L_{\overline{F}_j}(1-s),
\end{equation}
where
\[
L_{\overline{F}_j}(s):=\sum_{n=1}^{\infty}\frac{A_j(n,1)}{n^{s}},
\]
\begin{equation}\label{eq:Gfactors}
\begin{split}
G_{v,j}(s)
&:=\pi^{-\frac{3s}{2}}
\Gamma\!\left(\frac{s-\mu_{j,1}}{2}\right)
\Gamma\!\left(\frac{s-\mu_{j,2}}{2}\right)
\Gamma\!\left(\frac{s-\mu_{j,3}}{2}\right),
\\
\widetilde{G}_{v,j}(s)
&:=\pi^{-\frac{3s}{2}}
\Gamma\!\left(\frac{s+\mu_{j,1}}{2}\right)
\Gamma\!\left(\frac{s+\mu_{j,2}}{2}\right)
\Gamma\!\left(\frac{s+\mu_{j,3}}{2}\right).
\end{split}
\end{equation}
For details, see \cite[Chapter~6]{Go} or \cite{sun2024}.

\medskip
\noindent
For $s=\sigma+it$ with $t\neq 0$ and fixed $\sigma$, we have
\[
|G_{v,j}(s)| \asymp \prod_{l=1}^{3}
e^{-\left|\mathrm{Im}\left(it+\mu_{j,l}\right)\right|\pi/4}
\left|\mathrm{Im}\left(it+\mu_{j,l}\right)\right|^{\frac{\sigma-1}{2}}
\qquad
\text{(see \cite[(5.113)]{IK1}).}
\]
Note that when $|u| \ll |z|,$ 
\begin{equation}\label{gammadiff} 
\frac{d^{k}}{dz^{k}}\left(\frac{\Gamma(z+u)}{\Gamma(z)}\right) \asymp \left|\frac{\Gamma(z+u)}{\Gamma(z)} \left(\frac{u}{z}\right)^{k}\right|.
\end{equation}
\medskip
\noindent
We define the analytic conductor of $L_{F_j}(s)$ by
\[
q_{j,\infty}(\sigma+it):=\prod_{l=1}^{3}\bigl(3+|\sigma+it+\mu_{j,l}|\bigr).
\]
Let $q_j(s):=q_{j,\infty}(s)^{1/2}$ and $q_{j,\varepsilon}(s):=q_j(s)^{1+\varepsilon}$.
\subsection{The approximate functional equation}

\noindent
To apply the Kuznetsov formula, we use an approximate functional equation.

\begin{Lemma}\label{x5approximatefunctional}
 Let $s=\sigma+it$ with $\min(\mathrm{Re}(s),1-\mathrm{Re}(s))>0$, and let $T$ be sufficiently large. Then for any $t,t_1\in[T,2T]$,
\begin{equation}\label{eq:afe}
L_{F_j}(s)
=
\sum_{n=1}^{\infty}\frac{A_j(1,n)}{n^{s}}\,V_s\!\left(\frac{n}{q_{j,\varepsilon}(s)}\right)
+\frac{\widetilde{G}_{v,j}(1-s)}{G_{v,j}(s)}
\sum_{n=1}^{\infty}\frac{A_j(n,1)}{n^{1-s}}\,V'_{1-s}\!\bigl(q_{j,\varepsilon}(s)n\bigr),
\end{equation}
where
\begin{equation}\label{eq:Vsdefs}
\begin{split}
V_s(n)
&:=\frac{1}{2\pi i}\int_{(3)} n^{-u}\left(\cos\frac{\pi u}{4A}\right)^{-12A}\frac{G_{v,j}(s+u)}{G_{v,j}(s)}\frac{du}{u},
\\
V'_s(n)
&:=\frac{1}{2\pi i}\int_{(3)} n^{-u}\left(\cos\frac{\pi u}{4A}\right)^{-12A}\frac{\widetilde{G}_{v,j}(s+u)}{\widetilde{G}_{v,j}(s)}\frac{du}{u}.
\end{split}
\end{equation}
Moreover, for $\mathrm{Re}(s)>0$ we have
\begin{equation}\label{3V}
\begin{split}
V_s(n),\,V'_s(n) &\ll_{\mathcal{A}} \left(1+\frac{n}{q_j(s)}\right)^{-\mathcal{A}},\\
V_s(n),\,V'_s(n) &= 1+O\!\left(\left(\frac{n}{q_j(s)}\right)^{\frac{\mathrm{Re}(s)}{3}}\right),
\end{split}
\end{equation}
for any $\mathcal{A} >0.$
\end{Lemma}

\begin{proof}
See \cite[Theorem~5.3 and Proposition~5.4]{IK1}.
\end{proof}

\begin{Rem}
We choose the supports in \eqref{eq:afe} because the factor
\[
\frac{\widetilde{G}_{v,j}(1-s)}{G_{v,j}(s)}
\]
in the dual sum is not balanced with respect to the Langlands parameters $\mu_{j,l}$. This is the main reason that taking the two sums in the approximate functional equation to have comparable lengths does not work for the second moment of $L_{F_j}(s)M_{j;A,B}(s)$. Note also that directly applying an approximate functional equation with two dual sums of comparable lengths to $|L_{F_j}(s)|^{2}$ does not cause this issue; see \cite{sun2024}.
\end{Rem}

\subsection{Kuznetsov trace formula}
Given sufficiently large $U$, the test function $h_{U,M}$ satisfies the followings: 
\begin{enumerate}
\item
It holds that
\[
\sum_{\substack{F_j}}^{\ast} 1\asymp U^3M^{2}.
\]
See \cite[Corollary 7.2]{sun2024}.

\item (Variant of the Kuznetsov trace formula)
Suppose $|t|\le U^{1-\varepsilon}$. For any natural numbers $m_1,n_1,m_2,n_2$, and with $\delta(0)=1$ and $\delta(n)=0$ otherwise, we have
\begin{equation}\label{VKT}
\begin{split}
&\frac{\sum_{\substack{F_j}}^{\ast}
A_j(m_1,m_2)  A_j(n_2,n_1)}{H}
=\delta(n_1-m_1)\delta(n_2-m_2)
\\&+ O_{\varepsilon}\!\left(
(Un_1m_1n_2m_2)^{\varepsilon}
\Big(U^{-2}(m_1,m_2,n_1,n_2)^{7/64}+M^{-2}+U^{-2}M^{-2}(n_1m_1n_2m_2)^{1/2}\Big)
\right),
\end{split}
\end{equation}
(see \cite[Theorem 7.1]{sun2024}),
and when $n_{1} \neq m_{1}$ or $n_{2} \neq m_{2},$ we have 
\begin{equation}\label{VKT2}
\begin{split}
&\frac{\sum_{\substack{F_j}}^{\ast} 
A_j(m_1,m_2)A_j(n_2,n_1)V_{1/2+\theta+it}(m_{2})}{H} 
\\&= O_{\varepsilon}\!\left(
(Un_1m_1n_2m_2)^{\varepsilon}
\Big(U^{-2}(m_1,m_2,n_1,n_2)^{7/64}+M^{-2}+U^{-2}M^{-2}(n_1m_1n_2m_2)^{1/2}\Big)
\right).
\end{split}
\end{equation}

\end{enumerate}
\begin{Rem}
Attaching the weight function $V_{1/2+\theta+it}(m_2)$ does not change the size of the error term, since
\[
V_{1/2+\theta+it}(m_2) \ll 1,
\]
and the $k$th derivative (with respect to the Langlands parameter $\mu_{j,i}$) of
$V_{1/2+\theta+it}(m_2)$ satisfies
\[
\frac{\partial^k}{\partial \mu_{j,i}^k}
V_{1/2+\theta+it}(m_2)
\ll
\min_{l}\bigl(|\mu_{j,l}|\bigr)^{-k}
\]
(see \cite[Remark~6]{sun2024} and \cite{ButtBlom2020}).
This holds because the gamma factors
$G_{v,j}\!\left(\tfrac12+\theta+it+u\right)$ and
$G_{v,j}\!\left(\tfrac12+\theta+it\right)$
appearing in $V_{1/2+\theta+it}(m_2)$
have the same signs of $\mu_{j,i}$ for all $i$. See \eqref{gammadiff}.

\end{Rem}
\begin{Rem}
The best known bound for $|A_j(1,n)|$ is $|A_j(1,n)|\leq n^{5/14} d(n
)$. 
\end{Rem}
\section{Proof of Theorem \ref{Main}}

\begin{proof}
Let $s=\tfrac12+ \theta + it$. By Lemma~\ref{x5approximatefunctional}, for any $t\in[T,2T]$ we have
\begin{equation}\label{reductiontofinite}
L_{F_j}(s)
=
\sum_{m=1}^{\infty}\frac{A_j(1,m)}{m^{s}}
V_{s}\!\left(\frac{m}{q_{j,\varepsilon}(s)}\right)
+O\!\left(U^{-2022}\right).
\end{equation}
Therefore,
\begin{equation}\label{equation325}
\sum^{\ast}_{F_{j}}
\int_T^{2T}
L_{F_j}(1/2+\theta+it)\,M_{j;A,B}(1/2+\theta+it)\,dt
=
\sum^{\ast}_{F_{j}}\mathcal{S}_{j1}
+O\!\left(TU^{-1000}\right),
\end{equation}
where
\begin{equation}\label{eq:Sj1}
\begin{split}
\mathcal{S}_{j1}
:=
\sum_{m=1}^{\infty} m^{-1/2-\theta}
\sum_{\substack{n_{1}<A\\ n_{2}<B}}
\frac{\tau(n_{1})\tau(n_{2})}{n_{1}^{1/2+\theta}\,n_{2}^{\alpha/2+\alpha\theta}}
\int_T^{2T}
\,
\frac{A_j(1,m)\,A_j(n_2,n_1)}{(mn_1n_2^{\alpha})^{it}}
V_{1/2+\theta+it}\left(\frac{m}{q_{j,\varepsilon}(1/2+\theta+it)}\right)\,dt .
\end{split}
\end{equation}

\noindent
Here we frequently use the bound $q_j(1/2+\theta+it) \asymp U^{3/2}$.

\medskip
\noindent
To handle $\sum^{\ast}_{F_{j}}\mathcal{S}_{j1}$, we must consider both the diagonal contribution, corresponding to $n_2=m$ and $n_1=1$, and the off-diagonal contribution coming from the error term in the Kuznetsov trace formula. We show that the main term comes from the case $m=n_1=n_2=1$ in \eqref{eq:Sj1}, and that all other contributions are negligible.

\medskip
\noindent
We write
\[
\sum^{\ast}_{F_{j}}\mathcal{S}_{j1}=\mathcal{T}_D+\mathcal{T}_O,
\]
where $\mathcal{T}_D$ is the diagonal contribution and $\mathcal{T}_O$ is the off-diagonal contribution.
We first consider $\mathcal{T}_O$. We apply the Kuznetsov trace formula to
\[
\sum_{F_{j}}^{\ast} A_j(1,m)\,A_j(n_2,n_1)\,V_{1/2+\theta+it}(m).
\]
Truncating the sum over $m$ at $U^{3+6\varepsilon}$ (note that $q_j(s)q_{j,\varepsilon}(s)\ll U^{3+5\varepsilon}$), we obtain
\begin{equation}\label{eq:TO}
\begin{aligned}
\mathcal{T}_O
&\ll
U^{\varepsilon}
\sum_{m \le U^{3+6\varepsilon}} m^{-1/2}
\sum_{\substack{n_1 < A \\ n_2 < B}}
\frac{|\tau(n_1)\tau(n_2)|}{n_1^{1/2} n_2^{\alpha/2}}
\int_T^{2T}
\Bigl(
U(m n_1 n_2)^{1/2+\varepsilon}
+ U^{3}
+ U M^{2}(m n_1 n_2)^{7/64}
\Bigr)
\, dt
\\[0.4em]
&\ll
T U^{\varepsilon}
\sum_{\substack{n_1 < A \\ n_2 < B}}
|\tau(n_1)\tau(n_2)|\,
n_1^{-1/2} n_2^{-\alpha/2}
\Biggl[
(n_1 n_2)^{1/2+\varepsilon}
\sum_{m \le U^{3+6\varepsilon}} U m^{\varepsilon}
\\
&\hspace{3.2cm}
+ \sum_{m \le U^{3+6\varepsilon}} U^{3} m^{-1/2+\varepsilon}
+ \sum_{m \le U^{3+6\varepsilon}} U M^{2} m^{7/64 - 1/2 + \varepsilon}
(AB)^{7/64}
\Biggr]
\\[0.4em]
&\ll
T(\log U)^{O(1)}\, A^{1/2} B^{1-\alpha/2}
\Bigl[
(AB)^{1/2+\varepsilon}
U^{1+\varepsilon+(3+6\varepsilon)(1+\varepsilon)}
+ U^{\varepsilon + 3 + (3+6\varepsilon)(1/2+\varepsilon)}
\\
&\hspace{3.2cm}
+ (AB)^{7/64}
U^{1+\varepsilon+(3+6\varepsilon)(1/2+7/64+\varepsilon)} M^{2}
\Bigr]
\\[0.4em]
&\ll
T(\log U)^{O(1)}\Bigl(
A B^{3/2-\alpha/2} U^{4+20\varepsilon}
+ A^{1/2+\varepsilon} B^{1-\alpha/2} U^{9/2+10\varepsilon}
\\
&\hspace{2.8cm}
+ A^{7/64+1/2} B^{7/64+1-\alpha/2}
U^{1+117/64+20\varepsilon} M^{2}
\Bigr).
\end{aligned}
\end{equation}

\medskip
\noindent
Now we consider $\mathcal{T}_D$. Since $n_2<B$, the diagonal condition forces $m=n_2\le B$. Note that
\begin{equation}\label{asyv}
V_{1/2+\theta+it}(m)=1+O\!\left(\left(\frac{m}{q_j(s)}\right)^{1/6+\theta/3}\right).
\end{equation}
For the case $n_1=n_2=m=1$, the contribution equals
\[
\sum_{F_{j}}^{\ast}\int_T^{2T}
\left(1+O\!\left(\left(\frac{1}{q_j(s)q_{j,\varepsilon}(s)}\right)^{1/6+\theta/3}\right)\right)\,dt
=
\sum_{F_{j}}^{\ast} T\left(1+O(U^{-c})\right)
\]
for some constant $c>0$.

\medskip
\noindent
Next consider the case $n_1=1$ and $m=n_2$ with $n_2\ne 1$. Note that 
$$q_{j,\varepsilon}(1/2+\theta+it) q_{j}(1/2+\theta+it) \gg U^{3}.$$
For the delta term, we again use 
\eqref{asyv}. Therefore, the total contribution is
\begin{equation}\begin{split}
&\sum_{\substack{ n_{2}<B}}
\frac{\tau(n_{2})}{n_{2}^{\alpha/2+1/2+\theta+\alpha \theta}}
\int_T^{2T} \sum^{\ast}_{F_{j}}
\,
\frac{A_j(1,n_{2})\,A_j(n_2,1)}{(n_2^{1+\alpha})^{it}}
V_{1/2+\theta+it}\left(\frac{n_{2}}{q_{j,\varepsilon}(1/2+\theta+it)}\right)\,dt 
\\&=\quad \quad \sum_{\substack{ n_{2}<B}}
\frac{\tau(n_{2})}{n_{2}^{\alpha/2+1/2+\theta+\alpha\theta}}
\int_T^{2T} \sum^{\ast}_{F_{j}}
\,
\frac{A_j(1,n_{2})\,A_j(n_2,1)}{(n_2^{1+\alpha})^{it}}
\left(1+O\!\left(\left(\frac{n_{2}}{q_j(s)q_{j,\varepsilon}(s)}\right)^{1/6+\theta/3}\right)\right)\,dt 
\\& =\quad \quad \sum_{\substack{ n_{2}<B}}
\frac{\tau(n_{2})}{n_{2}^{\alpha/2+1/2+\theta+\alpha\theta}}
\int_T^{2T} \sum^{\ast}_{F_{j}}
\,
\frac{A_j(1,n_{2})\,A_j(n_2,1)}{(n_2^{1+\alpha})^{it}} dt 
\\& \quad \quad \quad +
O\left(  U^{-1/2-\theta}\sum_{\substack{ n_{2}<B}}
\frac{(\log U)^{O(1)}}{n_{2}^{\alpha/2+1/2+\theta+\alpha\theta}}
\int_T^{2T} \sum^{\ast}_{F_{j}}
n_{2}^{5/7+1/6+\theta/3} dt\right). \end{split}\end{equation}
By using the elementary bound 
\begin{equation}\label{eq:intbound}
\int_T^{2T} (n_2^{1+\alpha})^{-it}\,dt
\ll \left|\frac{1}{\log(n_2^{1+\alpha})}\right|,
\end{equation}
The delta term contribution in the first summation is bounded by 

\begin{equation}\label{eq:diagrest}
\sum^{\ast}_{F_{j}}
\sum_{n_{2}\le B}
\frac{(\log U)^{O(1)}}{n_{2}^{\alpha/2+1/2+\theta+\alpha\theta}|\log n_{2}|}
\ll(\log U)^{O(1)}
\sum^{\ast}_{F_{j}} \max\left(1,B^{1/2-\alpha/2 - \theta -\alpha \theta }\right).
\end{equation}
The non-delta term can be handled in the same way as for $\mathcal{T}_O$.
Now, let us consider the second summation
$$O\left(  U^{-1/2-\theta}\sum_{\substack{ n_{2}<B}}
\frac{(\log U)^{O(1)}}{n_{2}^{\alpha/2+1/2+\theta+\alpha\theta}}
\int_T^{2T} \sum^{\ast}_{F_{j}}
n_{2}^{5/7+1/6+\theta/3} dt\right),$$
which is bounded by 
\begin{equation}\begin{split} 
&U^{-1/2-\theta}T\sum^{\ast}_{F_{j}} \sum_{\substack{ n_{2}<B}}
\frac{(\log U)^{O(1)}}{n_{2}^{\alpha/2+1/2+\theta+\alpha \theta}}
n_{2}^{5/7+1/6+\theta/3} 
\\& \ll (\log U)^{O(1)}\sum^{\ast}_{F_{j}}U^{-1/2-\theta}T B^{1+5/7+1/6+\theta/3-\alpha/2-1/2-\theta-\alpha \theta}
\\& \ll(\log U)^{O(1)} \sum^{\ast}_{F_{j}} U^{-1/2-\theta}TB^{58/42-\alpha/2-2\theta/3-\alpha \theta}.
\end{split}
\end{equation}
Therefore, when $B^{58/42-\alpha/2-2\theta/3-\alpha \theta} \ll U^{1/2+\theta-\beta},$ the contribution is bounded by $T(\log U)^{O(1)}\sum^{\ast}_{F_{j}}U^{-\beta}.$ 
\end{proof}

\section{Proof of Theorem \ref{last}}
Note that in Theorem~\ref{Main}, the parameters $A$ and $B$ can be chosen arbitrarily small. Therefore, it suffices to consider only the square moment case. For this purpose, we use the following lemmas.

\begin{Lemma}\label{withoutmoment}
Let $0<\theta \leqslant 1/2$, $|t|<U^{1/4}$, and let $\ell_1,\ell_2<U^{1/3}$ be squarefree integers. Then, for any $\varepsilon>0$, we have
\[
\begin{aligned}
& \sum_{F_j}^{\ast} A_j(\ell_2,\ell_1)
\left|L_{F_j}\!\left(\tfrac12+\theta+it\right)\right|^2 \\
= {} &
\zeta(1+2\theta)\,
\frac{1}{(\ell_1\ell_2)^{1/2+\theta}}
\left(\frac{\ell_1}{\ell_2}\right)^{it} H \\
& + \zeta(1-2\theta)\,
\frac{1}{(\ell_1\ell_2)^{1/2-\theta}}
\left(\frac{\ell_1}{\ell_2}\right)^{it}
\frac{1}{192\pi^5}
\int_{\Re(\mu)=0}
h_{U,M}(\mu)
\prod_{k=1}^3
\left(-\frac{\mu_k^2}{4\pi^2}\right)^{-\theta}
\,\mathrm{d}_{\mathrm{spec}}\mu \\
& + O\!\left(
(\ell_1\ell_2)^{-1/2+\theta} U^{9/4-6\theta} M^2
+ U^{2-3\theta+\varepsilon} M^2
+ U^{83/42-41\theta/21+\varepsilon}
M^2(\ell_1\ell_2)^{7/64}
\right).
\end{aligned}
\]
\end{Lemma}

\begin{proof}
See \cite[Theorem~1.2]{sun2024}.
\end{proof}

Since $h_{U,M}$ is nonnegative and it is essentially supported for the case when $|\mu_{j,k}| \asymp |v_{j,k}| \asymp U$, we have
\[
\int_{\Re(\mu)=0}
h_{U,M}(\mu)
\prod_{k=1}^3
\left(-\frac{\mu_k^2}{4\pi^2}\right)^{-\theta}
\,\mathrm{d}_{\mathrm{spec}}\mu
\ll H U^{-6\theta}.
\]

Because we will take $\ell_1,\ell_2$ to be smaller than any fixed power of $U$, the error term will be negligible for our purposes.

We now expand $|M_{j;A,2}(s)|^{2}$ that is suitable for applying Lemma \ref{withoutmoment}.

\begin{Lemma}\label{mollifierexpansion}
Let
\[
M_{j;A,2}\!\left(\tfrac12+\theta+it\right)
= \sum_{n_1<A}
\frac{A_j(1,n_1)\tau(n_1)\mu(n_1)}
{n_1^{1/2+\theta+it}},
\]
for any function $\tau(n)$ that satisfies the conditions in Theorem \ref{Main}.
Then
\[
\begin{aligned}
\left|M_{j;A,2}\!\left(\tfrac12+\theta+it\right)\right|^2
&= \sum_{d<A} \frac{1}{d^{1+2\theta}}
\sum_{k_1<A/d}\sum_{k_2<A/d}
\frac{
A_j(k_2,k_1)\,
\tau(dk_1)\overline{\tau(dk_2)}\,
\mu(dk_1)\mu(dk_2)
}{
(k_1k_2)^{1/2+\theta}
}
\left(\frac{k_1}{k_2}\right)^{-it}.
\end{aligned}
\]
\end{Lemma}

\begin{proof}
By squaring out, 
\[
\left|M_{j;A,2}\!\left(\tfrac12+\theta+it\right)\right|^2
= \sum_{n_1<A}\sum_{n_2<A}
\frac{
A_j(1,n_1)A_j(n_2,1)\,
\tau(n_1)\overline{\tau(n_2)}\,
\mu(n_1)\mu(n_2)
}{
(n_1n_2)^{1/2+\theta}
}
\left(\frac{n_1}{n_2}\right)^{-it}.
\]
Using the Hecke relation
\[
A_j(1,n_1)A_j(n_2,1)
= \sum_{d\mid(n_1,n_2)} A_j\!\left(\tfrac{n_2}{d},\tfrac{n_1}{d}\right),
\]
and changing variables $n_1=dk_1$, $n_2=dk_2$ completes the proof.
\end{proof}

We are ready to prove the folloiwing theorem.
\begin{Th}\label{ending}
Let $\eta(U)$ be a function such that $(\log\log U)^2=o(\eta(U))$, and suppose that
\[
|t|=U^{\theta_1}, \qquad 
\frac{\eta(U)}{\log U}<\theta\leqslant \tfrac12,
\qquad 
A=\exp\!\left(\frac{\log U}{\eta(U)^{1/2}}\right),
\]
for some sufficiently small positive constant $\theta_1$.  
Let $M_{j,A,2}(1/2+\theta+it)$ be the Dirichlet polynomial defined in Lemma~\ref{mollifierexpansion}.  
Then
\[
\begin{aligned}
\sum_{F_j}^{\ast}
\left|
M_{j;A,2}\!\left(\tfrac12+\theta+it\right)
L_{F_j}\!\left(\tfrac12+\theta+it\right)
\right|^2
=
\zeta(1+2\theta)H
\sum_{d<A}\frac{1}{d^{1+2\theta}}
\sum_{k_1<A/d}\sum_{k_2<A/d}
\frac{
\tau(dk_1)\overline{\tau(dk_2)}
\mu(dk_1)\mu(dk_2)
}{
(k_1k_2)^{1+2\theta}
}
+o(H).
\end{aligned}
\]
\end{Th}

\begin{proof}
Substituting the identity in Lemma~\ref{mollifierexpansion} into Lemma~\ref{withoutmoment}, the first contribution is
\[
\zeta(1+2\theta)H
\sum_{d<A}\frac{1}{d^{1+2\theta}}
\sum_{k_1<A/d}\sum_{k_2<A/d}
\frac{
\tau(dk_1)\overline{\tau(dk_2)}
\mu(dk_1)\mu(dk_2)
}{
(k_1k_2)^{1+2\theta}
}.
\]
Since \[
\zeta(1-2\theta)\ll \frac{1}{2\theta},
\]
the second contribution is bounded by
\[
\begin{aligned}
\ll {} \theta^{-1}&
\sum_{d<A}\frac{1}{d^{1+2\theta}}
\sum_{k_1<A/d}\sum_{k_2<A/d}
\frac{
d(dk_1)d(dk_2)\,
|\mu(dk_1)\mu(dk_2)|
}{
k_1k_2
}
H U^{-6\theta}
\\
\ll \theta^{-1}{} &
H U^{-6\theta}(\log U)^{O(1)}.
\end{aligned}
\]
By the assumption $(\log\log U)^{2}=o(\eta(U))$, this is $o(H)$.

Since
\[
A=\exp\!\left(\frac{\log U}{\eta(U)^{1/2}}\right),
\]
the summation over $d,k_1,k_2$ are all smaller than any fixed power of $U$.  
Therefore, the contribution from the error term in Lemma~\ref{withoutmoment} is negligible.
\end{proof}
Note that the contribution from the term $n_{1}=1$ in \eqref{mollifierexpansion} is
$\zeta(1+2\theta)H.$
Now let us take $\tau(n)=C$ for some constant $C$ when $n \neq 1.$
Then the main term in Theorem~\ref{ending} becomes
\[  \zeta(1+2\theta)H + 
C^2 \zeta(1+2\theta)H\left(
\sum_{1<d<A}\frac{1}{d^{1+2\theta}}
\sum_{k_1<A/d}\sum_{k_2<A/d}
\frac{
\mu(dk_1)\mu(dk_2)
}{
(k_1k_2)^{1+2\theta}
}+\sum_{d<A}\frac{1}{d^{1+2\theta}}
\sum_{1<k_1<A/d}\sum_{1<k_2<A/d}
\frac{
\mu(dk_1)\mu(dk_2)
}{
(k_1k_2)^{1+2\theta}
}\right). 
\]
Using the trivial bound
\[
\sum_{k<m}\frac{\mu(m)}{m^{1+2\theta}}
=
\sum_{m=1}^{\infty}\frac{\mu(m)}{m^{1+2\theta}}
+O\!\left(\frac{m^{-2\theta}}{2\theta}\right),
\]
together with
\[
\zeta(1+2\theta)\ll \frac{1}{2\theta},
\]
 the main term is \begin{equation} \begin{split}    \zeta(1+2\theta)H +C^{2} & \zeta(1+2\theta)H\sum_{1<d<A} \frac{1}{d^{1+2\theta}}\left(\sum_{k=1}^{\infty} \frac{\mu(dk)}{k^{1+2\theta}} + O\left(\left(\frac{A}{d}\right)^{-2\theta}(2\theta)^{-1}\right)\right)^{2} 
 \\&+C^{2} \zeta(1+2\theta)H\sum_{d<A} \frac{1}{d^{1+2\theta}}\left(\sum_{k=2}^{\infty} \frac{\mu(dk)}{k^{1+2\theta}} + O\left(\left(\frac{A}{d}\right)^{-2\theta}(2\theta)^{-1}\right)\right)^{2}
 \\&=   \zeta(1+2\theta)H  +  C^{2}  \zeta(1+2\theta)H\sum_{1<d<A} \frac{1}{d^{1+2\theta}}\left(\sum_{k=1}^{\infty} \frac{\mu(dk)}{k^{1+2\theta}}\right)^{2}
 \\& + C^{2}  \zeta(1+2\theta)H\sum_{d<A} \frac{1}{d^{1+2\theta}}\left(\sum_{k=2}^{\infty} \frac{\mu(dk)}{k^{1+2\theta}}\right)^{2}+O\left(HA^{-2\theta} (\log U)^{O(1)}(2\theta)^{-3}\right)\end{split}\end{equation}

Since
\[
A^{-2\theta}(\log U)^{O(1)}(2\theta)^{-3}
=
\exp\!\left(
-2\theta\log A
+O(\log\log U)
-3\log(2\theta)
\right),
\]
the assumption on $\eta(U)$ implies that this error term is $o(H)$.
Therefore, by taking $C$ as 
\[
C=
\left(\frac{1-\zeta(1+2\theta)}{\zeta(1+2\theta)}
\left(\sum_{d=1}^{\infty}\frac{1}{d^{1+2\theta}}
\left(
\sum_{k=2}^{\infty}\frac{\mu(dk)}{k^{1+2\theta}}
\right)^{2}+\sum_{d=2}^{\infty}\frac{1}{d^{1+2\theta}}
\left(
\sum_{k=1}^{\infty}\frac{\mu(dk)}{k^{1+2\theta}}
\right)^{2} \right)
\right)^{-1/2},
\]
the main term is equal to $H$.

\bibliographystyle{plain}
\bibliography{over1}
\end{document}